\RequirePackage{fix-cm}
\documentclass[smallextended]{svjour3}       
\smartqed  
\usepackage{amssymb,epsfig,multirow}
\usepackage{graphicx}
\usepackage{epstopdf}
\usepackage{float}
\usepackage{color}
\usepackage{enumerate}
\usepackage{hyperref}
\usepackage[numbers,sort&compress]{natbib}
\usepackage{mathrsfs}
\usepackage{amsfonts}
\usepackage{amsmath}
\usepackage[noend]{algpseudocode}
\usepackage{algorithmicx,algorithm}

\newtheorem{prop}{Proposition}
\newtheorem{lem}{Lemma}
\newtheorem{ass}{Assumption}

\begin{document}
	
\title{Simultaneous perturbation stochastic approximation: towards one-measurement per iteration
\thanks{This research was supported by the Beijing Natural Science Foundation, grant Z180005 and by the National Natural Science Foundation of China under grants 12171021, and 11822103.}
}

\titlerunning{}        

\author{Shiru Li \and Yong  Xia \and Zi Xu
}


\institute{ S. Li \and    Y.  Xia\    \at
	 LMIB of the Ministry of Education, School of Mathematical Sciences, Beihang University, Beijing, 100191, P. R. China
	\email{lishiru@buaa.edu.cn  (S. Li); yxia@buaa.edu.cn (Y. Xia)}
	\and Z. Xu  \at
	Department of Mathematics, Shanghai University, Shanghai, 200444,  P. R. China
	\email{xuzi@shu.edu.cn (Z. Xu, corresponding author)}	
}

\date{Received: date / Accepted: date}

\maketitle

\begin{abstract}
When measuring the value of a function to be minimized is not only expensive but also with noise, the popular simultaneous perturbation stochastic approximation (SPSA) algorithm requires only two function values in each iteration. In this paper, we propose a method requiring only one function measurement value per iteration in the average sense. We prove the strong convergence and asymptotic normality of the new algorithm. Experimental results show the effectiveness and potential of our algorithm.


\end{abstract}

\section{Introduction}
We consider the following unconstrained optimization:
\begin{eqnarray*}
({\rm P})~~\min_{x\in\mathbb {R}^{n}}~f(x),
\end{eqnarray*}
where $f(x):\mathbb {R}^{n}\to\mathbb {R}$ is differentiable with the noisy measurements and $g(x)$ is the true gradient. In order to iteratively solve (P),
stochastic approximation (SA) is a popular algorithm scheme
given by
\begin{eqnarray*}
\hat x_{k+1}=\hat x_k-a_k d_k,
\end{eqnarray*}
where $\hat x_k$ is an estimation of the feasible solution at $k$-th iteration, $a_k\in\mathbb {R}$ is a step size and $d_k\in\mathbb {R}^n$ is a iterative direction.

If the gradient of the function is noisily available, by setting  $d_k$  as the noisy measurement of $g(x_k)$, which is given by
\begin{eqnarray*}
	\widetilde{g}_k=g(x_k)+\varepsilon_k,
\end{eqnarray*}
where $\varepsilon_k$ is the noise of the gradient of the $k$-th iteration,
SA  reduces to the Robbins-Monro (RM) algorithm \cite{RM1951SA}.

When the gradient is not available, the corresponding choice of the direction $d_k$ becomes an approximation of the gradient. Kiefer and Wolfowitz \cite{KW1952} proposed the finite difference stochastic approximation (FDSA) algorithm (also known as KW algorithm), which approximates $g(x)$ with the finite difference form $\hat g(x)$. That is, the $i$-th component of $\hat g(x_k)$ is given by
\begin{eqnarray*}
	\hat g_i(x_k)= \frac{\widetilde f(x_k+c_ke_i)-\widetilde f(x_k-c_ke_i)}{2c_k},
\end{eqnarray*}
where $\widetilde f(\cdot)$ is the measurement of $f(\cdot)$ with noise, $e_i$ is $i$-th column of the identity matrix and $c_k$ is a positive scalar.
FDSA algorithm needs $2n$ measurements of the function value in order to approximate a gradient vector.  Kushner and Clark \cite{1978Stochastic} proposed
the random direction stochastic approximation (RDSA) algorithm, which only requires two measurements to approximate $g(x)$:
\begin{eqnarray*}
	\hat g(x_k)= \frac{\widetilde f(x_k+c_k\xi_k)-\widetilde f(x_k-c_k\xi_k)}{2c_k}\xi_k,
\end{eqnarray*}
where $\xi_k$ is a random vector satisfying some specific distribution. Spall \cite{spall1992}  proposed the
simultaneous perturbation stochastic approximation (SPSA) method based on the following approximation
\begin{eqnarray}
	\hat g(x_k)= \frac{\widetilde f(x_k+c_k\xi_k)-\widetilde f(x_k-c_k\xi_k)}{2c_k}\xi^{-1}_k,\label{grad}
\end{eqnarray}
where $\xi^{-1}_k$ takes the inverse of every element of $\xi_k$,  $c_k$ is a positive scalar, and
\begin{eqnarray*}
	&\widetilde f(x_k+c_k\xi_k)&=f(x_k+c_k\xi_k)+\varepsilon_k^+,\\
	&\widetilde f(x_k-c_k\xi_k)&=f(x_k-c_k\xi_k)+\varepsilon_k^-,
\end{eqnarray*}
where $\varepsilon_k^+$ and $\varepsilon_k^-$ denote the measurement noise of the function value.
In \cite{1998Optimal}, Spall suggested an optimal choice of $\xi_k$ in SPSA by randomly, independently (also with respect to $\hat x_0,\hat x_1,\cdots,\hat x_k$) and uniformly generating in $\{-1,1\}^n$, i.e.,
the symmetric  Bernoulli distribution. The following settings of the stepsize $a_k$ and the perturbation parameter $c_k$
\begin{eqnarray}
	a_k&=&\frac{a}{(k+1+A)^\alpha}, \label{ak}\\
	c_k&=&\frac{c}{(k+1)^\gamma},  \label{ck}
\end{eqnarray}
are due to Spall \cite{1998Implementation}, where $a$, $A$, $c$, $\alpha$ and $\gamma$ are predefined constants.

Assume that $\varepsilon_k^{+},~\varepsilon_k^{-}$, satisfy
\begin{eqnarray*}
	&~&E\left[ \varepsilon_k^{+ }-\varepsilon_k^{- }\big|\mathscr{F}_k,\xi_k\right] =0~a.s.,~\forall k,\\
	&~&\mathscr{F}_k\equiv\{\hat x_0,\hat x_1,\cdots,\hat x_k\},
\end{eqnarray*}
where $E[\cdot]$ stands for the expectation, and a.s. represents almost surely.
It can be proved that $\hat g_k(x)$ is an unbiased estimation of $g(x)$ so that $\hat g_k(x)$ can be regarded as a good approximation of $g(x)$. Under the above assumptions,  strong convergence and asymptotic normality of the iterations $\hat x_k$ for SPSA have been established  in \cite{spall1992}.

SPSA only needs {\it two} measurements of objective function values. For  problems  of $n$ dimension, the number of functional measurements in each iteration of SPSA is $n$ times less than that of FDSA. This superiority makes SPSA very popular,  with widespread applications in  control engineering, signal processing, neural network training, parameter estimation, etc.. Many variants and improvements of SPSA are developed, for example, the second-order SPSA \cite{spall1995second}, the accelerated SPSA \cite{spall1997accelerated,zhu2002modified}, SPSA for nonsmooth optimization  \cite{2007nonsmooth}, the adaptive direction version \cite{xu2008adaptive}, and the fuzzy adaptive SPSA \cite{2011fuzzy}.

In order to further improve SPSA from two to one functional measurement per iteration, Spall \cite{spall1997one} presented a one-measurement version,
where
\begin{equation*}
	\hat g(x_k)= \frac{\widetilde f(x_k+c_k\xi_k)}{2c_k}\xi^{-1}_k.
\end{equation*}
This algorithm, denoted by SPSA1, was reported to could outperform the classical SPSA in some special cases. Although convergence and asymptotic normality results of  SPSA1 have been established, there is a bias term in the asymptotic covariance matrix of  SPSA1 in the iterative convergence process. It makes the practical performance of  SPSA1 not as good as expected. Consequently, the following fundamental problem remains open: \\

{\it Is there an efficient SPSA algorithm with only one function value measurement in each iteration?} \\

The difficulty is that with one measurement of function value one can not  approximate the gradient properly. In this paper, we present a new algorithm, which
evaluates two function values for every two iterations. We established its strong convergence and asymptotic normality. Numerical experiments demonstrate the efficiency comparing with SPSA and SPSA1.
Therefore, our algorithm can be regarded as an efficient SPSA with only one functional measurement per iteration {\it in the average sense}.

The rest of this paper is as follows. We propose a new SPSA algorithm, and then establish its  strong convergence and asymptotic normality in Section 2. Section 3 reports numerical results demonstrating the efficiency of our algorithm.  Conclusions are made in Section 4.

{\bf Notation.} Denote by $f^*$ the minimum value of the function $f$.
$\hat g_k(x_k)$ is often simply rewritten as $\hat g_k$. Let $H(x)$ be the Hessian matrix for $f$. Denote by $\|\cdot\|$ the Euclidean norm. The tensor product is denoted by $\otimes$.
$O(\cdot)$ stands for the infinitesimal of the same order of $(\cdot)$.
$C_{n}^{m}~(n,m\in \mathbb{R}$, and $m\le n)$ stands for the  combination number given by
$$C_{n}^{m}=\frac{n(n-1)\cdots (n-m+1)}{m(m-1)\cdots 1}.$$
Let $\Omega=\{\omega\}$ be the sequence $\hat x_1,\hat x_2,\cdots$ generated by the sample space. $I$ denotes the identity matrix.
$\text{Diag}(\cdot)$ represents the diagonal matrix with $(\cdot)$ being the diagonal elements. $\text{sgn}(\cdot)$ returns the element-wise sign vector of $(\cdot)$. $\mathscr{I}_{\{\cdot\}}$ denotes the indicator function. Weak convergence (convergence in distribution) is denoted by $\stackrel{w}{\longrightarrow}$. We abbreviate ``with probability one'' to w.p.1..

\section{One-measurement SPSA algorithm}
In this section, we positively answer the question raised above by presenting a new version of SPSA  algorithm with one measurement per iteration in the sense of average. Theoretical analysis is also provided.

\subsection{The new algorithm}

Our idea is based on three observations for SPSA algorithm:
\begin{itemize}
\item[(a)] The search direction $\xi_k$ in $k$-th iteration is independent of any information of the objective function $f(\cdot)$.
\item[(b)] The difference of two functional measurements $f(x_k\pm c_k\xi_k)$ decides which side of the search direction $\xi_k$ is descent.
\item[(c)] The stepsizes of the adjacent two iterations is close to each other.
\end{itemize}
Therefore, a good prediction of the descent side of $\xi_k$ can replace the two functional measurements in one iteration. Our motivation is the following fact:

{\it Suppose $-\hat g_k$ is the steepest descent, then $-\hat\xi_k$ is a direction of descent if $\hat g^T_k\hat\xi_k\geq0$.}

  Our algorithm first takes one step along the direction $-\hat g_k$ and then one step along $-\hat\xi_k$. Two function measurements are required in the first step, and no function measurement is required in the second step. In this way, only two function measurement points are needed for every two steps, which is equivalent to only one function measurement value for one step in the average sense. Based on the above description, we named the method SPSA1-A. Algorithm \ref{alg1} introduces the framework of our method.
\begin{algorithm}[htbp]
	\caption{SPSA1-A}\label{alg1}
	\hspace*{0.02in} {\bf Input:}  \hspace*{0.02in}
	initial values $\hat x_0,~a_0,~c_0$\\
	\hspace*{0.02in} {\bf Output:} $\hat x_{k+1}$
	\begin{algorithmic}[1]
		\State Compute $$\rho= \begin{cases}
		 C_{n-1}^{n/2}/\left(2^{n-1}+ C_{n}^{{n}/{2}}/2\right),~ &\text{if}~n~\text{is even},  \\
			 C_{n-1}^{{(n-1)}/{2}}/2^{n-1},~ &\text{otherwise}.\\
		\end{cases}$$
		\For{$k=0,1,\cdots$}
		\State Randomly and uniformly generate $\xi_k\in\{-1,1\}^n$.
		\State Calculate  $a_k$ and $c_k$ (for example,  \eqref{ak}-\eqref{ck}).
Compute  $\widetilde{f}(x_k\pm c_k\xi_k)$, $\hat g_k$ by \eqref{grad}, and $\rho_k=\rho/\|\hat g_k\|_{\infty}$.
		\State Update $\hat x_{k+\frac{1}{2}}\leftarrow\hat x_k-a_k\hat g_k/(1+\rho_k).$
		\State Randomly and uniformly generate $\hat\xi_k\in\{d\in\{-1,1\}^n:~ d^{T}\hat g_k\ge 0\}$.
		\State Update $\hat x_{k+1}\leftarrow\hat x_{k+\frac{1}{2}}- a_k\hat\xi_k/(1+\rho_k)$.
		\State Stop if the termination criterion is reached.
		\EndFor
		\State \Return $\hat x_{k+1}$
	\end{algorithmic}
\end{algorithm}
The termination criterion in Algorithm SPSA1-A can be employed as either the maximum number of iterations or the solution accuracy.

\subsection{Theoretical analysis}
In this section we establish the strong convergence and the asymptotic normality for   Algorithm SPSA1-A. We leave the proofs to Appendix as they are similar to that in \cite{1998Optimal}.

Lemma \ref{lem1} reveals the relationship between $\hat\xi_k$ and $\hat g_k$ in the sense of conditional mathematical expectations.
\begin{lem}\label{lem1}
	\begin{eqnarray*}
		E[ \hat{\xi}_k|\hat{x}_k]  = E\left[  \rho_k\hat g_k |\hat x_k\right].
	\end{eqnarray*}
\end{lem}

Lemma \ref{lem2} shows that the bias of the estimate of $(\hat g_k(\cdot)+\hat\xi_k)/(1+\rho_k)$ as $g(\cdot)$ goes to $0$ as $k\rightarrow\infty$ under certain conditions.

\begin{lem}\label{lem2}
Suppose there is an index $K<\infty$ and a constant $r>0$ such that $f^{(3)}(x)\triangleq\partial^3f/(\partial x^T)^3$ is not only continuous but also element-wise bounded in $\{x:~\|x-\hat x_k\|<r\}$ for all $k\geq K$.
Then for almost all $\omega\in\Omega$, we have
	\begin{eqnarray*}
		b_k(\hat x_k)&\triangleq&E\left[  ( \hat g_k+\hat\xi_k ) /(1+\rho_k)-g(\hat x_k)\Big|\hat x_k\right]
 =  O(c_k^2) ~(c_k\rightarrow 0).
	\end{eqnarray*}
\end{lem}

In order to complete the proof of convergence, we define the solution error
\begin{eqnarray*}
	e_k(\hat x_k) = (\hat g_k + \hat\xi_k )/(1+\rho_k)- E\left[ (\hat g_k + \hat\xi_k)/(1+\rho_k)\Big|\hat x_k\right],
\end{eqnarray*}
and then we have
\begin{eqnarray*}
	\hat x_{k+1}=\hat x_{k}-a_k\left[ g_k(\hat x_k)+b_k(\hat x_k)+e_k(\hat x_k)\right] \label{GRM}.
\end{eqnarray*}

We now present some necessary assumptions.

\begin{ass}\label{ass1}
	$a_k,~c_k > 0~\forall~k;~a_k\to0,~c_k\to0~as~k\to\infty;$
	\begin{eqnarray*}
		\sum_{k = 0}^{\infty}a_k = \infty,~\sum_{k = 0}^{\infty}\bigg( \frac{a_k}{c_k}\bigg) ^2< \infty.
	\end{eqnarray*}
\end{ass}
Assumption \ref{ass1} allows more choices of $a_k$ and $c_k$ rather than \eqref{ak}-\eqref{ck}.

\begin{ass} \label{ass2}
There exist $\alpha_0>0,~\alpha_1> 0$ such that $\forall ~k,~E\left(\varepsilon_k^{\pm}\right)^2 \le \alpha_0,$ $E\left(f{(\hat x_k \pm c_k\xi_k)}\right)^2 \le \alpha_1$.
\end{ass}

\begin{ass}
	$\sup_k\|\hat x_k\| < \infty ~a.s.$.
\end{ass}

\begin{ass}\label{ass4}
	$x^*$ is an asymptotically stable solution of the differential equation $d\phi(t)/dt = -g(x)$.
\end{ass}

\begin{ass}\label{ass5}
	Let $D(x^*) = \{\phi_0:\lim\limits_{t\to\infty}\phi(t|\phi_0) = \theta^*\}$
	where $\phi(t|\phi_0)$ denotes the solution to the differential equation of Assumption \ref{ass4} based on initial conditions $\phi_0$ (i.e., $D(x^*)$ is the domain of attraction). There exists a compact $S\subseteq D(x^*)$ such that $\hat x_k\in S$ infinitely often for almost all sample points.
\end{ass}

\begin{prop}\label{prop1}
Under Assumptions \ref{ass1}-\ref{ass5} and conditions of Lemma \ref{lem2}, when $k\to\infty$, it holds that
	\begin{eqnarray*}
		\hat x_k\to x^*~~~~for~almost~all~\omega\in\Omega.
	\end{eqnarray*}
\end{prop}

Similar to \cite{spall1992}, we can now establish the asymptotic normality analysis for Algorithm SPSA1-A.
For the sake of simplicity, we let $a_k=a/k^\alpha$ and $c_k=c/k^\gamma$ where $a,c,\alpha,\gamma> 0$. We strengthen Assumption \ref{ass2} as the following one:
\begin{ass} \label{ass2'}
There exist $\delta>0,~\alpha_0>0,~\alpha_1>0$ such that $\forall~k,~E\left| \varepsilon^{\pm}_k\right| ^{2+\delta}\leq\alpha_0,~E|f(\hat x_k\pm c_k\xi_k)|^{2+\delta}\leq\alpha_1$.
\end{ass}


\begin{prop}\label{prop2}
Suppose the assumptions made in Lemma \ref{lem2} and Proposition \ref{prop1} hold,  and  Assumption \ref{ass2} holds as a replacement of Assumption \ref{ass2'}. Let $\sigma^2$ be such that $E\left[ ( \varepsilon_k^{+}-\varepsilon_k^{-}) ^2\big|\mathscr{F}_k\right] \rightarrow\sigma^2$ a.s., as $k\rightarrow\infty$. For sufficiently large $k$ and almost all $\omega $, let the sequence $E\left[ (\varepsilon_k^{+}-\varepsilon_k^{-})^2\big|\mathscr{F}_k, c_k\xi_k=\eta\right] $ be equicontinuous at $\eta=0$ and continuous with respect to $\eta$ in a compact, connected set containing $c_k\xi_k$ a.s.. Furthermore, let $\beta=\alpha-2\gamma>0,~3\gamma-\alpha/2\geq0$, and $P$ be orthogonal with $PH(x^*)P^{T}=a^{-1}diag(\lambda_1,\cdots,\lambda_n)$. Then, we have
	\begin{eqnarray*}
		k^{\beta/2}(\hat x_k-x^*)\stackrel{w}{\longrightarrow}\mathscr{N}
		,~k\rightarrow\infty,
	\end{eqnarray*}
	where $\mathscr{N}$ is a Gaussian random vector with $E\mathscr{N}=\mu$ and $Cov(\mathscr{N},\mathscr{N})=PMP^T$,
	$M=\frac{1}{4} a^2c^{-2}\sigma^2diag\left[ (2\lambda_1-\beta_+)^{-1},\cdots,(2\lambda_n-\beta_+)^{-1}\right] $ with $\beta_+=\beta<2\min_i\lambda_i$ if $\alpha=1$ and $\beta_+=0$ if $\alpha<1$,
	\begin{eqnarray*}
		\mu= \begin{cases}
			0,~&\text{if}~3\gamma-\alpha/2>0, \\
			(aH(x^*)-\frac{1}{2}\beta_+I)^{-1}T,~&\text{if}~3\gamma-\alpha/2=0, \\
		\end{cases}
	\end{eqnarray*}
and	the $l$-th element of T is given by
	\begin{eqnarray*}
		-\frac{1}{6}ac^2\Bigg[ f^{(3)}_{lll}(x^*)+3\sum_{\substack{i=1\\i\neq l}}^{n}f_{iil}^{(3)}(x^*)\Bigg] .
	\end{eqnarray*}
\end{prop}

\section{Numerical experiments}
In this section we do numerical experiments to compare Algorithms SPSA1-A  with  SPSA. We also numerically compare Algorithms  SPSA1-A and SPSA1 due to Spall \cite{spall1997one}. For all the tested algorithms, we start from the same initial point $x_0$ and stop when the termination criterion is reached. The optimal values of all test  functions are all zero. We sample the noise $\varepsilon_k$ from a normal distribution with mean $0$ and standard variance $0.01$.

We set $a_k$ and $c_k$ in Algorithm SPSA as \eqref{ak}-\eqref{ck}.
When Algorithm SPSA1-A stops in finite steps, $1+\rho_k$ is upper bounded w.p.1.. All assumptions for theoretical analysis are satisfied if we set $\{a_k\}$ to $\{a_k(1+\rho_k)\}$. So in practical version of Algorithm SPSA1-A, we set
\begin{eqnarray*}
	a_k&=&\frac{a(1+\rho_k)}{(k+1+A)^\alpha}.
\end{eqnarray*}

\subsection{Test $\uppercase\expandafter{\romannumeral1}$}
We first test the following three unconstrained minimization problems given in \cite{1981Testing}.
For each example, we empirically choose the parameters for Algorithm SPSA to be the best, 
and then independently run each algorithm $50$ times to get an average iterative curve.
We terminate algorithms when the maximum number
$M$ of iterations is reached.

Problem $1$. (Rosenbrock function)
\begin{equation*}
	f(x)=100(x_2-x_1^2)^2+(1-x_1)^2
\end{equation*}
We choose $x_0=(-1.2,1)^{T}$ and $M=5000$. The parameters given in Table \ref{tab1} are empirically the best for Algorithm SPSA.
\begin{table}[htbp]\centering
	\caption{The parameters in Algorithms SPSA and SPSA1-A for Rosenbrock function.} \label{tab1}
	\begin{center}
		\begin{tabular}{cccccc} \hline
			$$  &a  &A  &c &$\alpha$ &$\gamma$    \\
			\hline
			SPSA & 0.1& 2200& 0.1 &0.602 &0.101\\
			SPSA1{\rm-}A &0.1 &2200 &0.1 &0.602 &0.101\\			
			\hline			
		\end{tabular}
	\end{center}
\end{table}
\begin{figure}[htbp]
	\centering
	\includegraphics[width=6.66cm,height=4.86cm]{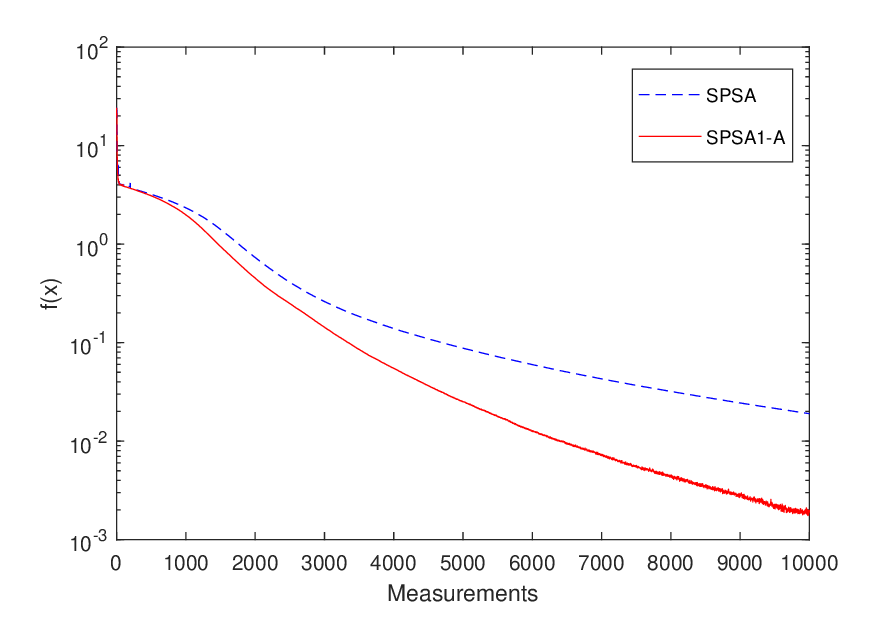}
	\caption{Minimizing the Rosenbrock function by Algorithms SPSA and SPSA1-A.}\label{fig1}
\end{figure}

As shown in Figure \ref{fig1}, the logarithmic plot of the average number of function measurements, Algorithm SPSA1-A  converges faster than SPSA. In particular, to output a solution of the same accuracy, say $f(x)<0.01$, Algorithm SPSA1-A requires less than half as many functional measurements as that of SPSA.

Problem $2$. (Beale function)
\begin{equation*}
	f(x)=[1.5-x_1(1-x_2)]^2+[2.25-x_1(1-x_2^2)]^2+[2.625-x_1(1-x_2^3)]^2
\end{equation*}
We choose $x_0=(1,1)^{T}$ and $M=5000$. The parameters given in Table \ref{tab2} are empirically the best for Algorithm SPSA.
\begin{table}[htbp]\centering
	\caption{The parameters in Algorithms SPSA and SPSA1-A for Beale function.} \label{tab2}
	\begin{center}
		\begin{tabular}{cccccc} \hline
			$$  &a  &A  &c &$\alpha$ &$\gamma$     \\
			\hline
			SPSA &1 &30 &0.1 &1 &0.16667\\
			SPSA1{\rm-}A &1 &30 &0.1 &1 &0.16667\\	
			\hline			
		\end{tabular}
	\end{center}
\end{table}
\begin{figure}[htbp]
	\centering
	\includegraphics[width=6.66cm,height=4.86cm]{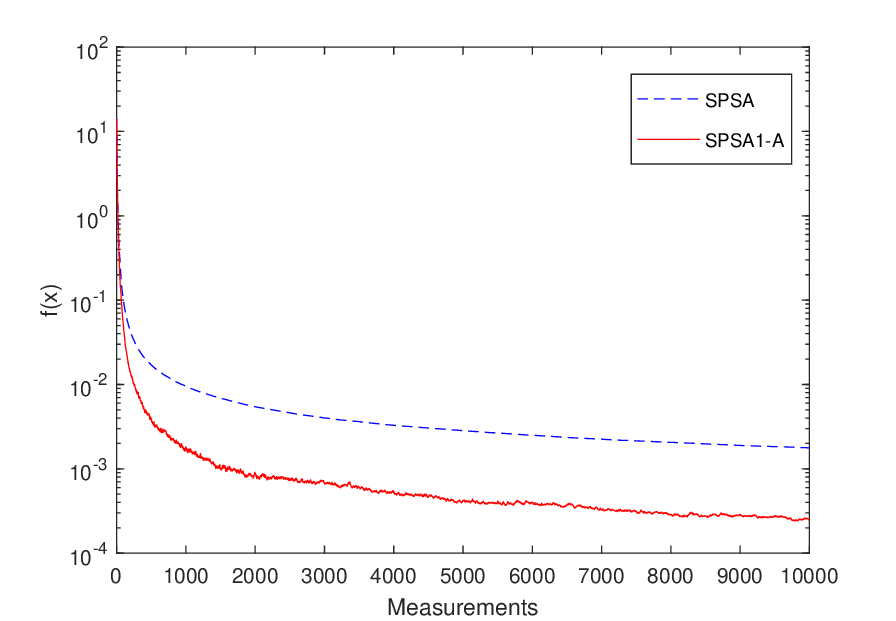}
	\caption{Minimizing the Beal function by Algorithms SPSA and SPSA1-A.}\label{fig2}
\end{figure}

According to Figure \ref{fig2}, the logarithmic plot of the average number of function measurements,  the iterative functional values by Algorithm SPSA1-A  decreased much faster than that of SPSA in the first hundreds of functional measurements. After $M$ iterations, the  accuracy of the  solution outputted  by Algorithm SPSA1-A is an order of magnitude higher than  that of Algorithm SPSA.

Problem $3$. (Powell singular function)
\begin{equation*}
	f(x)=(x_1+10x_2)^2+5(x_3-x_4)^2+(x_2-2x_3)^4+10(x_1-x_4)^4
\end{equation*}
We choose $x_0=(3,-1,0,1)^{T}$ and $M=5000$. The parameters given in Table \ref{tab3} are empirically the best for Algorithms SPSA and SPSA1-A, respectively.
\begin{table}[htbp]\centering
	\caption{The parameters in Algorithms SPSA and SPSA1-A for Powell singular function.} \label{tab3}
	\begin{center}
		\begin{tabular}{cccccc} \hline
			$$  &a  &A  &c &$\alpha$ &$\gamma$     \\
			\hline
			SPSA &0.08 &1000 &0.1 &0.602 &0.101\\
			SPSA1{\rm-}A &0.02 &100 &0.1 &0.602 &0.101\\	
			\hline			
		\end{tabular}
	\end{center}
\end{table}

\begin{figure}[htbp]
	\centering
	\includegraphics[width=6.66cm,height=4.86cm]{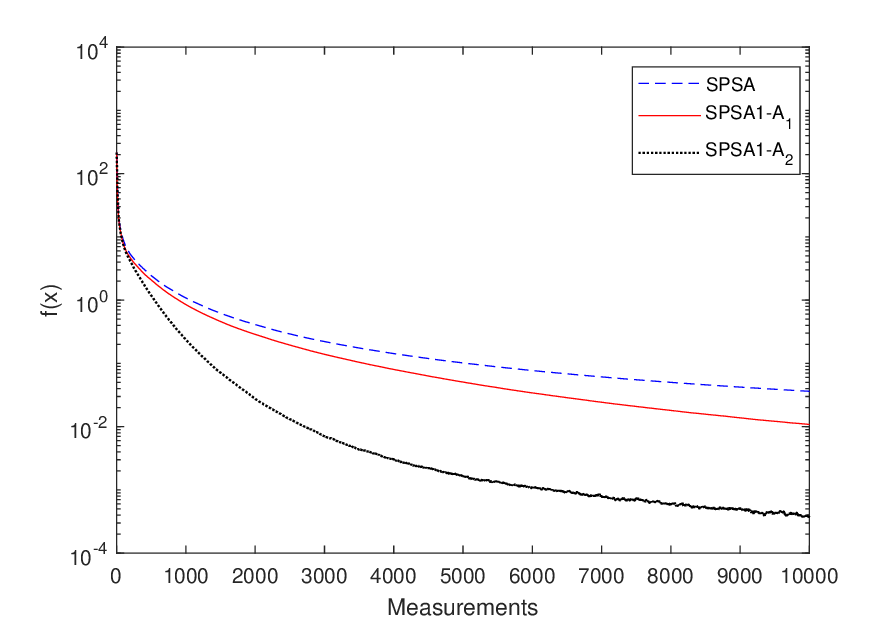}
	\caption{Minimizing the Powell singular function by Algorithms SPSA and SPSA1-A.}\label{fig3}
\end{figure}

Algorithms SPSA1-A$_1$ and SPSA1-A$_2$ use the best parameters for Algorithms SPSA and SPSA1-A, respectively. As shown in the logarithm plot \ref{fig3}, both algorithms converge faster than Algorithm  SPSA. Moreover, as expected, Algorithm SPSA1-A$_2$ highly outperforms SPSA1-A$_1$.

\subsection{Test $\uppercase\expandafter{\romannumeral2}$}
We numerically compare Algorithms SPSA, SPSA1-A and SPSA1 on the example presented for  SPSA1 in \cite{spall1997one}:
\begin{equation}\label{case2}
	f(x)=x^Tx+0.1\sum_{i=1}^{5}x_i^3+0.01\sum_{i=1}^{5}x_i^4.
\end{equation}
All algorithms start from the same initial point $x_0=(3,-1,0,1)^{T}$ and stop when either the maximum iteration number $M=10^5$ is reached or the error $e_k=|f(x_k)-f^*|$ is less than a threshold (see the last two columns of
Tables \ref{tab4} and \ref{tab5}). We test all algorithms with the same two kinds of parameters as that in \cite{spall1997one}, see Column 2-6 in Tables \ref{tab4} and \ref{tab5}.
\begin{table}[htbp]\centering
	\caption{Numerical results for minimizing (\ref{case2}) based on the parameters of the first kind. } \label{tab4}
	\begin{center}
		\begin{tabular}{cccccccc} \hline
			$$  &a  &A  &c &$\alpha$ &$\gamma$ &$e_k\le 10^{-2}$ &$e_k\le 10^{-3}$  \\
			\hline
			SPSA &0.17 &20 &0.06 &1 &0.16667 &206 &7711\\
			SPSA1 &0.17 &20 &0.06 &1 &0.16667 &3930 &--\\
			SPSA1{\rm-}A &0.17 &20 &0.06 &1 &0.16667 &80 &784\\
			\hline			
		\end{tabular}
	\end{center}
\end{table}

\begin{table}[htbp]\centering
	\caption{Numerical results for minimizing (\ref{case2}) based on the parameters of the second kind.} \label{tab5}
	\begin{center}
		\begin{tabular}{cccccccc} \hline
			$$  &a  &A  &c &$\alpha$ &$\gamma$ &$e_k\le 10^{-2}$ &$e_k\le 10^{-3}$  \\
			\hline	
			SPSA &0.27 &100 &0.06 &1 &0.16667 &349 &3738\\
			SPSA1 &0.27 &100 &0.06 &1 &0.16667 &2172 &--\\
			SPSA1{\rm-}A &0.27 &100 &0.06 &1 &0.16667 &144 &711\\
			\hline			
		\end{tabular}
	\end{center}
\end{table}

We report the numbers of functional measurements of three algorithms until they terminate in the last two columns of Tables \ref{tab4} and  \ref{tab5}, where ``--'' stands for the situation that the maximum iteration number $M$ is reached.
Clearly, Algorithm SPSA1-A performs much better than the other two algorithms. We also notice that Algorithm SPSA1 fails to find solution of high accuracy.

\section{Conclusion}
The simultaneous perturbation stochastic approximation (SPSA) algorithm is  popular  for minimizing a noised function. It measures two function values in each iteration. It makes sense that each iteration requires at least two functional measurements to guarantee the descent of the search direction. In this paper, we propose a new algorithm measuring two function values every two iterations, that is, only one measurement of function is taken per iteration in the average sense. We prove the strong convergence and asymptotic normality of the new algorithm. Numerical results demonstrate the effectiveness of our new algorithm comparing with Algorithm SPSA. Future works include more applications and further improvement of our new algorithm.

\bibliographystyle{plain}

\appendix
\section{Proof of  Lemma \ref{lem1}.}
\begin{proof}
We start from the observation
			$$E[ \hat\xi_k|\hat x_k]  = E [ E[\hat\xi_k|\hat x_k,\hat g_k] |  \hat x_k].$$
Notice that each component $\hat\xi_{ki}~(i=1,2,\cdots,n)$ of $\hat\xi_k$ obeys the symmetric Bernoulli distribution and satisfies ${\hat{\xi}_k}^T \hat g_k\ge0$. Therefore, at least half of the components of $\hat\xi_k$ have the same sign as the components $\hat g_{ki}~(i=1,2,\cdots,n)$ of $\hat g_k$. If $n$ is even, $\hat\xi_k$ has $C_{n}^{ n/2}+C_{n}^{ n/2 +1}+\cdots+C_{n}^{n}=2^{n-1}+ C_{n}^{n/2}/2$ choices. So for any possible choice $\zeta_j$, we have  \[
P(\hat\xi_k=\zeta_j)=\frac{1}{2^{n-1}+ C_{n}^{ n/2}/2}.
\]
	
If the signs of $\hat\xi_{ki}$ and $\hat g_{ki}$ $(\forall i\in \{1,2,\cdots,n\})$ are the same, then at least $ n/2-1$ of the remaining $n-1$  elements of $\hat\xi_k$ share the same signs as that of $\hat g_{ki}$. In this case, $\hat\xi_{k}$ has $C_{n-1}^{(n-2)/2}+C_{n-1}^{ n/2}+\cdots+C_{n-1}^{n-1}=2^{n-2}+C_{n-1}^{(n-2)/2}$ choices. Then we can write
	\begin{eqnarray*}
		\Sigma\zeta_{j}&=&\left\lbrace 2^{n-2}+C_{n-1}^{(n-2)/2}-\left[ \left( 2^{n-1}+ C_{n}^{n/2}/2\right) -\left( 2^{n-2}+C_{n-1}^{(n-2)/2}\right) \right] \right\rbrace \text{sgn}\left( \hat g_{k}\right) \\
		&=&\left( 2C_{n-1}^{(n-2)/2}- C_{n}^{n/2}/2\right) \text{sgn}\left( \hat g_{k}\right) \\
		&=&C_{n-1}^{n/2}\text{sgn}\left( \hat g_{k}\right).
	\end{eqnarray*}
	
	If $n$ is odd, $\{\zeta_j\}$ has $C_{n}^{(n+1)/2}+C_{n}^{(n+1)/2+1}+\cdots+C_{n}^{n}=2^{n-1}$ choices, each with probability $P(\hat\xi_k=\zeta_j)= 1/2^{n-1}$. The sign of $\hat\xi_{ki}$ is either the same as or opposite to
that of	$\hat g_{ki}$. If their signs are the same, then at least $(n-1)/2$ of the remaining $n-1$  elements of
	$\hat\xi_k$ share the same signs as $\hat g_{ki}$. In this case, $\hat\xi_k$ has
	$C_{n-1}^{(n-1)/2}+C_{n-1}^{(n+1)/2}+\cdots+C_{n-1}^{n-1}=2^{n-2}+ C_{n-1}^{(n-1)/2}/2$ choices.
	Then we can write
	\begin{eqnarray*}
		\Sigma\zeta_{j}&=&\left\lbrace 2^{n-2}+ C_{n-1}^{(n-1)/2}/2-\left[ 2^{n-1}-\left( 2^{n-2}+ C_{n-1}^{(n-1)/2}/2\right) \right] \right\rbrace \text{sgn}(\hat g_{k})\\
		&=&C_{n-1}^{(n-1)/2}\text{sgn}\left( \hat g_{k}\right).
	\end{eqnarray*}
It then holds that
	$$\Sigma\zeta_{j}P(\hat\xi_k=\zeta_j)= \begin{cases}
		\frac{C_{n-1}^{n/2}}{2^{n-1}+ C_{n}^{n/2}/2}\text{sgn}(\hat g_{k}),~ &{\rm if}~n~{\rm is~even},    \\
		\frac{C_{n-1}^{(n-1)/2}}{2^{n-1}}\text{sgn}(\hat g_{k}),~ &{\rm otherwise}.\\
	\end{cases}$$
For $\rho$ defined in Step $1$ of Algorithm SPSA1-A,
we have
		\begin{eqnarray*}
			&~&E [E[\hat\xi_k|\hat x_k,\hat g_k] | \hat x_k]
			 = E [_{\hat g_k}\rho E[\text{sgn}(\hat g_k)|\hat x_k] | \hat x_k ]
			 =\rho E[\text{sgn}(\hat g_k)|\hat x_k].
	\end{eqnarray*}
Let $\rho_k=\rho/\|\hat g_k\|_{\infty}$. We can complete the proof
	$$E[\hat\xi_k|\hat x_k]=E[\rho_k\hat g_k|\hat x_k].$$
\end{proof}


\section{ Proof of Lemma \ref{lem2}.}
\begin{proof}
By a proof similar to that of Lemma~\ref{lem1}, we have
	\begin{eqnarray*}
		E\left[ \frac{\hat\xi_k}{1+\rho_k}\bigg|\hat x_k\right]  =E\left[ \frac{\rho_k}{1+\rho_k}\hat g_k\bigg|\hat x_k\right].
	\end{eqnarray*}
	Then, we can obtain
	\begin{eqnarray*}
		b_k(\hat x_k)&=& {E[\hat g_k(\hat x_k) - g(\hat x_k)|\hat x_k]}.
	\end{eqnarray*}
In the sequel,  the lemma can be proved similar to that in \cite{spall1992}.
\end{proof}

\section{Proof of Proposition \ref{prop1}. }

\begin{proof}
	According to \cite{1978Stochastic}, based on Lemma \ref{lem2} and Assumption \ref{ass1}, we have

	i)
	\begin{eqnarray*}
		\|b_k(\hat x_k)\|<\infty ~\forall~k~{\rm and}~b_k(\hat x_k)\to 0,~a.s..
	\end{eqnarray*}
Next, we prove that
	
	ii)
	\begin{eqnarray*}
		\lim\limits_{k\to\infty}P\left( \sup\limits_{m\ge k}\left\| \sum_{i=k}^{m}a_ie_i(\hat x_i)\right\| \ge\eta\right) =0,~ {\rm for~any}~\eta>0.
	\end{eqnarray*}
First, for any $l\in\{1,2,\cdots,n\}$, we have
	\begin{eqnarray*}
		Ee_{il}^{2}&\le& Var\left( \frac{\hat g_{il}+\hat\xi_{il}}{1+\rho_i}\right)
		\le E(\frac{\hat g_{il}+\hat\xi_{il}}{1+\rho_i})^2\\
		&=&E\frac{\hat g_{il}^2}{(1+\rho_i)^2}+2E\frac{\hat g_{il}\hat\xi_{il}}{(1+\rho_i)^2}
		+E\frac{\hat\xi_{il}^2}{(1+\rho_i)^2}.
	\end{eqnarray*}
	Then we obtain
	\begin{eqnarray*}
		E\left( \hat g_{il}\hat\xi_{il}\right) &\le&E\left( \left| \hat g_{il}\right| \cdot\left| \hat\xi_{il}\right| \right) \\
	& \le & E \left| \hat g_{il}\right|  \\
		&\le&\frac{1}{2c_i}E\left[ \left| \xi_{il}^{-1}\right| \cdot
		\left( \left| L(\hat x_i+c_i\xi_i)\right| +\left| L(\hat x_i-c_i\xi_i)\right| +\left| \varepsilon_i^+\right| +\left| \varepsilon_i^-\right| \right) \right] \\
		&\le&\left( \sqrt{\alpha_0}+\sqrt{\alpha_1}\right) c_i^{-1}.
	\end{eqnarray*}
It follows from \cite{spall1992} that
	\begin{eqnarray*}
		E\hat g_{il}^2\le2(\alpha_1+\alpha_0)c_i^{-2}.
	\end{eqnarray*}
Since $\frac{1}{(1+\rho_i)^2}=\frac{1}{(1+\rho/|g_{il}|)^2}\le 1$ and
	$E\hat\xi_{il}^2= 1$, we have
	\begin{eqnarray*}
		Ee_{il}^2&\le& E\hat g_{il}^2+2E\hat g_{il}\hat\xi_{il}+E\hat\xi_{il}^2
		 \le 2(\alpha_1+\alpha_0)c_i^{-2}+(\sqrt{\alpha_0}+\sqrt{\alpha_1})c_i^{-1}+1.
	\end{eqnarray*}
Therefore, it holds that
\[
E\|e_k\|^2\le p[2(\alpha_1+\alpha_0)c_k^{-2}+(\sqrt{\alpha_0}+\sqrt{\alpha_1})c_k^{-1}+1].
\]
	
Next, since $\{\sum_{i=k}^{m}a_ie_i\}_{m\ge k}$ is a martingale sequence, it follows from the inequality in \cite[P. 315]{1953Stochastic} (see also \cite[P. 27]{1978Stochastic}) that
\begin{eqnarray}
		P(\sup\limits_{m\ge k}\|\sum_{i=k}^{m}a_ie_i\|\ge\eta\|)&\le&\eta^{-2}E\|\sum_{i=k}^{\infty}a_ie_i\|^2
		 = \eta^{-2}\sum_{i=k}^{\infty}a_{i}^{2}E\|e_i\|^2, \label{p2}
\end{eqnarray}
where the equality holds as $E\left[ e_{i}^Te_j\right] =E\left[ e_{i}^TE\left[ e_j|\hat x_j\right] \right] =0,~\forall~i<j$.

Then, by (\ref{p2}) and Assumption \ref{ass1}, we complete the proof of ii).
\end{proof}

\section{Proof of Proposition \ref{prop2}.}

\begin{proof}
	In order to complete the proof, we need to verify whether  conditions (2.2.1), (2.2.2), and (2.2.3) in Fabian \cite{fabian1968} are true. Here we assume that all assumptions on $\theta_k$ or $\mathscr{F}_k$ hold.
	According to the notation in \cite{fabian1968}, we can get
	\begin{eqnarray*}
		\hat x_{k+1}- x^*=(I-k^{-\alpha}\Gamma_k)(\hat x_k-x^*)+k^{-(\alpha+\beta)/2}\Phi_k V_k+k^{-\alpha-\beta/2}T_k,
	\end{eqnarray*}
	where $\Gamma_k=aH~(\overline{x}_k)$, $V_k=k^{-\gamma}\left\lbrace \frac{\hat g_k(\hat x_k)+\hat\xi_k}{1+\rho_k}-E\left[ \frac{\hat g_k(\hat x_k)+\hat\xi_k}{1+\rho_k}\Big|\hat x_k\right] \right\rbrace $, $\Phi_k=-aI$, and $T_k=-ak^{\beta/2}b_k(\hat x_k)$.
	In fact, there is an open neighborhood of $\hat x_k$ (for $k$ sufficiently large) containing $ x^*$  in which $H(\cdot)$ is continuous. Then
	\begin{eqnarray*}
		E\left[ \frac{\hat g_k(\hat x_k)+\hat\xi_k}{1+\rho_k}\bigg|\hat x_k\right] &=&H(\overline{x}_k)(\hat x_k-x^*)+b_k(x^*),
	\end{eqnarray*}
	where $\Gamma_k=aH(\overline{x}_k)$ lies in the line segment between $\hat x_k$ and $x^*$.
	
Based on the continuity of $H(\cdot)$ and a.s. convergence of $\hat x_k$, we have $\Gamma_k=aH(\overline{x}_k)\rightarrow aH(x^*)$ a.s..
	
	Now we prove the convergence of $T_k$ for $3\gamma-\alpha/2\ge 0$. When $3\gamma-\alpha/2>0$, as $b_k(\hat x_k)=O(k^{-2\gamma})$ a.s., we can write that $T_k\rightarrow0$ a.s.. When $3\gamma-\alpha/2=0$, by the facts that $\hat x_k\rightarrow x^*$ a.s. and the uniformly boundedness of $L^{(3)}$ near $x^*$,
we have
	\begin{eqnarray*}
		k^{2\gamma}b_{kl}(\hat x_k)-\frac{1}{6}c^2L^{(3)}(x^*)E(\xi_k\otimes\xi_k\otimes\xi_k)\rightarrow 0~a.s..
	\end{eqnarray*}
	
	Then ${\xi_{ki}}$ is symmetrically i.i.d. for each $k$, which means that the $l$-th element of $T_k$ satisfies that
	\begin{eqnarray*}
T_{kl}\rightarrow-\frac{1}{6}ac^2{L_{lll}^{(3)}(x^*)+\sum_{\substack{i=1\\i\neq l}}^{p}\left[ L^{(3)}_{lii}(x^*)+L_{ili}^{(3)}(x^*)+L_{iil}^{(3)}(x^*)\right] }~a.s..
	\end{eqnarray*}
Therefore,  $T_k$ converges for $3\gamma-\alpha/2\geq0$.
	
We can write
	\begin{eqnarray*}
 E\left[ V_kV_k^T\big|\mathscr{F}_k\right ]
		&=&k^{-2\gamma}\Bigg\lbrace E\left[ \frac{1}{(1+\rho_k)^2}(\hat g_k+\hat\xi_k)(\hat g_k+\hat\xi_k)^T\bigg|\hat x_k\right]	\\	
& & ~-E\left[ \frac{1}{1+\rho_k}(\hat g_k+\hat\xi_k)\bigg|\hat x_k\right] E\left[ \frac{1}{1+\rho_k}(\hat g_k+\hat\xi_k)^T\bigg|\hat x_k^T\right] \Bigg\rbrace,
	\end{eqnarray*}
	where
	\begin{eqnarray*}
		&~&k^{-2\gamma}E\bigg[ \frac{1}{(1+\rho_k)^2}(\hat g_k+\hat\xi_k)(\hat g_k+\hat\xi_k)^T\bigg|\hat x_k\bigg] \\
		&=&k^{-2\gamma}E\left[ \frac{1}{(1+\rho_k)^2}(\hat g_k\hat g_k^T+\hat g_k\hat\xi_k^T+\hat\xi_k\hat g_k^T+\hat\xi_k\hat\xi_k^T)\bigg|\hat x_k\right] .\\
	\end{eqnarray*}	
Define $\xi^{-1}_k:=(\xi^{-1}_{k1},\cdots,\xi^{-1}_{kp})^T$. Then we have
	\begin{eqnarray}\label{E[ggT]}
 E\left[ \frac{1}{(1+\rho_k)^2}\hat g_k\hat g_k^T\bigg|\mathscr{F}_k\right]
		&=&E\left\lbrace \frac{1}{(1+\rho_k)^2}\xi_k^{-1}(\xi_k^{-1})^T\left[ \frac{\varepsilon_k^{(+)}-\varepsilon_k^{(-)}}{2ck^{-\gamma}}\right] ^2\bigg|\mathscr{F}_k\right\rbrace,
	\end{eqnarray}
	where
	\begin{eqnarray*}
		\frac{1}{(1+\rho_k)^2} = \frac{1}{(1+\frac{\rho}{|\hat g_{kl}|})^2}
		 = \frac{|\hat g_{kl}|^2}{(|g_{kl}|+\rho)^2} = \frac{[\frac{y^{(+)}-y^{(-)}}{2ck^{-\gamma}}]^2}{[\frac{y^{(+)}-y^{(-)}}{2ck^{-\gamma}}+\rho]^2}
		 = \frac{[y^{(+)}-y^{(-)}]^2}{[y^{(+)}-y^{(-)}+2ck^{-\gamma}\rho]^2},
	\end{eqnarray*}
and	the last equation $\rightarrow1$ with $k\rightarrow\infty$. Therefore,  (\ref{E[ggT]}) is same as the third term in (3.5) in \cite{spall1992}. As the element of $\hat g_k\hat\xi_k^T+\hat\xi_k\hat g_k^T+\hat\xi_k\hat\xi_k^T$ is bounded,  we have
	\begin{eqnarray*}
		k^{-2\gamma}E\bigg[\frac{1}{(1+\rho_k)^2}\left( \hat g_k\hat\xi_k^T+\hat\xi_k\hat g_k^T+\hat\xi_k\hat\xi_k^T\right) \bigg|\hat x_k\bigg]\rightarrow0,
	\end{eqnarray*}
	and
	\begin{eqnarray*}
 E\left[ \frac{1}{1+\rho_k}(\hat g_k+\hat\xi_k)\bigg|\hat x_k\right] E\left[ \frac{1}{1+\rho_k}(\hat g_k+\hat\xi_k)^T\bigg|\hat x_k^T\right]
		 = E\left( \hat g_k\big|\hat x_k\right) E\left( \hat g_k^T\big|\hat x_k^T\right) .
	\end{eqnarray*}
	According to \cite{spall1992}, we obtain
	\begin{eqnarray*}
		E\left( V_kV_k^T\big|\mathscr{F}_k\right) \rightarrow\frac{1}{4}c^{-2}\sigma^2I~~~~a.s..
	\end{eqnarray*}
Thus we have obtained the conditions (2.2.1) and (2.2.2) of \cite{fabian1968}. Next we prove condition (2.2.3), i.e.,
	\begin{eqnarray*}
		\lim\limits_{k\rightarrow\infty}E\big(\mathscr{I}_{\{\left\| V_k\right\| ^2\geq rk^\alpha\}}\left\| V_k\right\| ^2\big)=0~~~\forall r>0.
	\end{eqnarray*}
	By Holder's inequality and $0<\delta'<\delta/2$, the upper bound of the above limit can be obtained as
	\begin{eqnarray*}
		&~&\lim\limits_{k\rightarrow\infty}\sup P(\left\| V_k\right\| ^2\geq rk^\alpha)^{\delta'/(1+\delta')}(E\left\| V_k\right\| ^{2(1+\delta')})^{1/(1+\delta')}\\
		&\leq&\lim\limits_{k\rightarrow\infty}\sup \left( \frac{E\left\| V_k\right\| ^2}{rk^\alpha}\right) ^{\delta'/(1+\delta')}(E\left\| V_k\right\| ^{2(1+\delta')})^{1/(1+\delta')}.
	\end{eqnarray*}
	Notice that
	\begin{eqnarray*}
		\left\| V_k\right\| ^{2(1+\delta')}\leq 2^{2(1+\delta')}k^{-2(1+\delta')\gamma}\bigg[ \left\| \frac{\hat g_k+\hat\xi_k}{1+\rho_k}\right\| ^2+\left\| b_k\right\| ^2+\left\| g_k\right\| ^2\bigg].
	\end{eqnarray*}
Then the proof is completed following from the proof of \cite[Proposition 1]{spall1992}.
\end{proof}

\end{document}